\documentclass[preprint]{elsarticle}

\usepackage{hyperref}

\journal{Discrete Applied Mathematics}

\usepackage[latin1]{inputenc}
\usepackage[T1]{fontenc}
\usepackage{ae}
\usepackage{latexsym,ifthen}
\usepackage{amsmath}
\usepackage{amssymb}
\usepackage{epsfig}
\usepackage{graphicx}
\usepackage{float}
\restylefloat{figure}
\usepackage{xspace}
\usepackage{tikz}

\newtheorem{Proposition}{Proposition}
\newtheorem{Corollary}{Corollary}
\newtheorem{Theorem}{Theorem}

\newcommand{\proof}{\noindent\textit{Proof.}}
\newcommand{\N}{\mathbb{N}}
\newcommand{\Probe}{\mathbb{P}}

\newcommand{\INPUT}{\textbf{Input: }}
\newcommand{\OUTPUT}{\textbf{Output: }}
\newcommand{\IF}{\textbf{if }}
\newcommand{\THEN}{\textbf{then }}
\newcommand{\ELSE}{\textbf{else }}
\newcommand{\ENDIF}{\textbf{endif}}

\newlength{\ai}
\newcommand{\ei}{\hspace*{\ai}}

\begin{document}

\begin{frontmatter}

\title{Good characterizations and linear time recognition for 2-probe block graphs\tnoteref{mytitlenote}}
\tnotetext[mytitlenote]{A preliminary version appeared in the proceedings of International Computer Symposium 2014, Frontiers in Artificial Intelligence and Applications 274, pp. 22--31, IOS Press, 2015.}


\author[mymainaddress]{Van Bang Le}
\ead{van-bang.le@uni-rostock.de}

\author[mysecondaryaddress]{Sheng-Lung Peng\corref{mycorrespondingauthor}}
\cortext[mycorrespondingauthor]{Corresponding author}
\ead{slpeng@mail.ndhu.edu.tw}

\address[mymainaddress]{Universit\"at Rostock, Institut f\"ur Informatik, Germany}
\address[mysecondaryaddress]{Department of Computer Science and Information Engineering,\\
  National Dong Hwa University,
  Shoufeng, Hualien 97401, Taiwan}

\begin{abstract}
Block graphs are graphs in which every block (biconnected component) is a clique. A graph $G=(V,E)$ is said to be an (unpartitioned) $k$-probe block graph if there exist $k$  independent sets $\N_i\subseteq V$, $1\le i\le k$, such that the graph $G'$ obtained from $G$ by adding certain edges between vertices inside the sets $\N_i$, $1\le i\le k$, is a block graph; if the independent sets $\N_i$ are given, $G$ is called a partitioned $k$-probe block graph.
In this paper we give good characterizations for $2$-probe block graphs, in both unpartitioned and partitioned cases. As an algorithmic implication, partitioned and unpartitioned probe block graphs can be recognized in linear time, improving a recognition algorithm of cubic time complexity previously obtained by Chang {\em et al.} [Block-graph width, Theoretical Computer Science 412 (2011), 2496--2502].
\end{abstract}

\begin{keyword}
Probe graph\sep block graph\sep probe block graph\sep probe complete graph
\end{keyword}

\end{frontmatter}


\section{Introduction}

Given a graph class $\mathcal{C}$, a graph $G=(V,E)$ is called a \emph{probe $\mathcal{C}$ graph} if there exists an independent set $\N\subseteq V$ (of \emph{non-probes}) and a set $E'\subseteq\binom{\N}{2}$ such that the graph $G'=(V,E\cup E')$ is in
the class $\mathcal{C}$, where $\binom{\N}{2}$ stands for the set of all 2-element subsets of $\N$.
A graph $G=(V,E)$ with a \emph{given} independent set $\N\subseteq V$ is said to be a \emph{partitioned probe $\mathcal{C}$ graph} if there exists a set $E'\subseteq\binom{\N}{2}$ such that the graph $G'=(V,E\cup E')$ is in the class $\mathcal{C}$. In both cases, $G'$ is called a $\mathcal{C}$ {\em embedding} of $G$. Thus, a graph is a (partitioned) probe $\mathcal{C}$ graph if and only if it admits a $\mathcal{C}$ embedding.

Recognizing partitioned probe $\mathcal{C}$ graphs is a special case of the $\mathcal{C}$-\textsc{graph sandwich} problem (cf.~\cite{GolKapSha95}). More precisely, given two graphs $G_i=(V, E_i)$, $i=1,2$, on the same vertex $V$ such that $E_1\subseteq E_2$, the $\mathcal{C}$-\textsc{graph sandwich} problem asks for the existence of a graph $G=(V,E)$ such that $E_1\subseteq E\subseteq E_2$ and $G$ is in $\mathcal{C}$. Recognizing partitioned probe $\mathcal{C}$ graphs with a given independent set $\N$ is a special case of the $\mathcal{C}$-\textsc{graph sandwich} problem, where $E_2\setminus E_1=\binom{\N}{2}$. Both concepts stem from computational biology; see, {\em e.g.\/},~\cite{GolKapSha94,GolKapSha95,Zhang,Zhang-etal}.

Probe graphs have been investigated for various graph classes; see~\cite{ChaChaKloPen} for more information.

Recently, the concept of probe graphs has been generalized in~\cite{ChaHunKloPen}.
A graph $G$ is said to be a \emph{$k$-probe $\mathcal{C}$ graph} if there exist independent sets $\N_1, \ldots, \N_k$ in $G$ such that there exists a graph $G'\in \mathcal{C}$ (an \emph{embedding} of $G$) such that for every edge $xy$ in $G'$ which is not an edge of $G$ there
exists an $i$ with $x, y \in \N_i$. In the case $k = 1$, $G$ is a \emph{probe $\mathcal{C}$ graph}.

We refer to the partitioned case of the problem when a collection of independent sets
$\N_i$, $i = 1,\ldots, k$, is a part of the input; otherwise, it is an unpartitioned case.
For historical reasons we call the set of vertices $\Probe = V \setminus \bigcup^k_{i=1} \N_i$ the set of \emph{probes} and the vertices of $\bigcup^k_{i=1} \N_i$ the set of
\emph{non-probes}.

In~\cite{ChaHunKloPen}, $k$-probe complete graphs and $k$-probe block graphs have been investigated. The authors proved that, for fixed $k$, $k$-probe complete graphs can be characterized by finitely many forbidden induced subgraphs, their proof is however not constructive. They also showed, implicitly, that $k$-probe complete graphs and $k$-probe block graphs can be recognized in cubic time. The case $k=1$, {\em e.g.}, probe complete graphs and probe block graphs, has been discussed in depth in~\cite{LePeng12}.

In this paper, we study $2$-probe complete graphs and $2$-probe block graphs in more details. Our main results are:
\begin{itemize}
 \item A characterization of partitioned $2$-probe block graphs in terms of certain ``enhanced graph'' (Theorem~\ref{thm:p2-probeblock}), stating that $G$ is a partitioned $2$-probe block graph if and only if the enhanced graph $G^*$ is a block graph.
 \item Forbidden induced subgraph characterizations of unpartitioned $2$-probe block graphs (Theorem~\ref{thm:2-probeblock}).
 \item Linear time recognition for $2$-probe block graphs, in both partitioned and unpartitioned cases.
\end{itemize}
The first result is of great interest because the enhanced graph contains only necessary edges, {\em i.e.}, new edges that must be added. In this sense, the enhanced graph is an optimal embedding of the probe graph.
This type of characterization is rarely possible, and our result is the first one in case $k=2$. In case of probe graphs, {\em i.e.}, $k=1$ only few are known: In~\cite{BayLedeR} it is shown that a graph is a partitioned probe threshold graph, respectively, a partitioned probe trivially perfect graph if and only if a certain enhanced graph is a threshold graph, respectively, a trivially perfect graph. In~\cite{Le} it is shown that a graph is a partitioned chain graph if and only if a certain enhanced graph is a chain graph, and recently, \cite{LePeng12} (cf. Theorem~\ref{thm:pprobeblock}) proved that a graph is a partitioned block graph if and only if a certain enhanced graph is a block graph. For some other cases, a certain enhanced graph can be defined that admits some nice properties; see~\cite{CohGolLipSte,GolLip,Zhang}.

Forbidden induced subgraph characterizations are very desirable as they (or their proofs)
often imply polynomial time for recognition, and give a lot of structural
information of the graphs.\footnote{That is why characterizing probe interval graphs by forbidden induced subgraphs is a long-standing interesting open problem; see~\cite{MWZ}} This is the case with the second result. Based on our forbidden induced subgraph characterization, we will obtain a linear time algorithm for recognizing if a given graph is a $2$-probe block graph, improving the cubic time complexity provided previously in~\cite{ChaHunKloPen}.

The paper is structured as follows. In Section~\ref{sec:def}, we collect all the necessary
definitions, and review results about probe complete graphs and probe block graphs. In Section~\ref{sec:2-probecomplete}, we discuss $2$-probe complete graphs. Partitioned and unpartitioned $2$-probe block graphs will be considered in Section~\ref{sec:p2-probeblock} and in Section~\ref{sec:2-probeblock}, respectively. A linear time recognition algorithm of unpartitioned 2-probe block graphs is proposed in Section~\ref{sec:reog}. We conclude the paper with some open problems in Section~\ref{sec:conclusion}.

\section{Definitions and notion}\label{sec:def}
In a graph, a set of vertices is an {\em independent set}, respectively, a {\em clique} if no two, respectively, every two vertices in this set are adjacent. For two graphs $G$ and $H$, we write $G+H$ for the disjoint union of $G$ and $H$, and $2G$ for $G+G$. The {\em join} $G\star H$ is obtained from $G+H$ by adding all possible edges $xy$ between any vertex $x$ in $G$ and any vertex $y$ in $H$. The complete graph with $n$ vertices is denoted by $K_n$. The path and cycle with $n$ vertices of length $n-1$, respectively, of length $n$, is denoted by $P_n$, respectively, $C_n$. Let $G = (V, E)$ be a graph. For a vertex $v\in V$ we write $N(v)$ for the set of its neighbors in $G$. A \emph{universal} vertex $v$ is one such that $N(v)\cup\{v\}=V$.
For a subset $U\subseteq V$ we write $G[U]$ for the subgraph of $G$ induced by $U$ and $G-U$ for the graph $G[V\setminus U]$; for a vertex $v$ we write $G-v$ rather than $G[V\setminus \{v\}]$.

A (connected or not) graph is a \emph{block graph} if each of its maximal 2-connected components, {\em i.e.\/}, its blocks, is a clique. A \emph{chordal graph} is one in which every cycle $C_\ell$ of length $\ell\ge 4$ has a chord. (A chord of a cycle is an edge not belonging to the cycle but joining to vertices of the cycle.) A \emph{diamond} is the complete graph on four vertices minus an edge. It is well-known (and easy to see) that block graphs are exactly the chordal graphs without induced diamond.

\begin{Proposition}[Folklore]
A graph is a block graph if and only if it is a diamond-free chordal graph.
\end{Proposition}

Here, given a graph $F$, a graph is said to be {\em $F$-free} if it has no induced subgraph isomorphic to $F$.  For a set of graphs $\mathcal{F}$, a graph is said to be {\em $\mathcal{F}$-free} if it is $F$-free for each $F\in \mathcal{F}$.

A graph $G$ is called \emph{distance-hereditary} if for all vertices $u,v \in V(G)$ any induced path between $u$ and $v$ is a shortest path. A graph $G$ is called \emph{ptolemaic} if, in any connected component of $G$, every four vertices 
 satisfy the so-called \emph{ptolemaic inequality} (cf.~\cite{Howorka}). 

\begin{Proposition}[Folklore]\label{dh}
\mbox{}
\begin{itemize}
\item[\em (i)] Ptolemaic graphs, gem-free chordal graphs, and $C_4$-free distance-hereditary graphs coincide.
\item[\em (ii)] Distance-hereditary graphs are exactly the graphs without induced house, hole, domino, gem.
\end{itemize}
\end{Proposition}

Here, a \emph{house} is a $5$-cycle with exactly one chord, a \emph{hole} is a $C_\ell$, $\ell\ge 5$, a \emph{domino} is a $6$-cycle with exactly one long chord, and a \emph{gem} is the join $P_4\star K_1$.

Another graph class that will be important in our discussion is the class of $P_4$-free graphs, or \emph{cographs}. Clearly, by Proposition~\ref{dh}, cographs are distance-hereditary, and we will often use the following well known fact.

\begin{Proposition}[Folklore]
Any connected cograph $G$ is the join $G=G_1\star G_2$ of two smaller cographs $G_1, G_2$.
\end{Proposition}

For graph classes not defined here see, for example,~\cite{BraLeSpi,ChaChaKloPen,Golumbic}.


A \emph{split graph} is a graph whose vertex set can be partitioned into a clique and an independent set. It is well-known that split graphs are exactly the chordal graphs without induced $2K_2$.

\begin{Proposition}[\cite{FolHam}]\label{prop:split}
A graph is a split graph if and only if it is a $2K_2$-free chordal graph.
\end{Proposition}

A \emph{complete split graph} is a split graph $G=(V,E)$ admitting a partition $V=Q\cup S$ into a clique $Q$ and an  independent set $S$ such that every vertex in $Q$ is adjacent to every vertex in $S$. Such a partition is also called a \emph{complete split partition} of a split graph. Note that if the complete split graph $G=(V,E)$ is not a clique, then $G$ has exactly one complete split partition $V=Q\cup S$.

\begin{Proposition}[\cite{LePeng12}]\label{prop:completesplit}
The following statements are equivalent for any graph $G$.
\begin{itemize}
 \item[\em (i)] $G$ is a probe complete graph;
 \item[\em (ii)] $G$ is a $\{K_2+K_1, C_4\}$-free graph;
 \item[\em (iii)] $G$ is a $(K_2+K_1)$-free split graph;
 \item[\em (iv)] $G$ is a complete split graph.
\end{itemize}
\end{Proposition}

Given a graph $G=(V,E)$ together with an independent set $\N\subseteq V$, the \emph{enhanced graph} $G^*=(V,E^*)$ is obtained from $G$ by adding all edges between two vertices in $\N$ that are two vertices of an induced diamond in $G$.

Partitioned probe block graphs can be characterized as follows.

\begin{Theorem}[\cite{LePeng12}]\label{thm:pprobeblock}
Let $G=(V,E)$ be a graph with a partition $V= \Probe\cup \N$, where $\N$ is an independent set.
Then the following statements are equivalent: 
\begin{itemize}
 \item[\em (i)] $G=(\Probe\cup \N, E)$ is a partitioned probe block graph;
 \item[\em (ii)] $G$ is a ptolemaic graph and satisfies the property that the two non-adjacent vertices of every induced diamond in $G$ belong to $\N$;
 \item[\em (iii)] Every block $B$ of $G$ is a complete split graph with $B=(B\cap \Probe) \cup (B\cap \N)$ a complete split partition;
 \item[\em (iv)] $G^*$ is a block graph.
\end{itemize}
\end{Theorem}

Note that condition (ii) in Theorem~\ref{thm:pprobeblock} above can be equivalently stated using three partitioned induced forbidden diamonds (cf.~\cite{LePeng12}).

Probe block graphs can be characterized as follows; see Fig.~\ref{fig:F1234} for the graphs $F_1, F_2$ and $F_3$.
\begin{Theorem}[\cite{LePeng12}]\label{thm:probeblock}
The following statements are equivalent for any graph $G$:
\begin{itemize}
 \item[\em (i)] $G$ is a probe block graph;
 \item[\em (ii)] $G$ is an $\{F_1, F_2$, $F_3\}$-free ptolemaic graph;
 \item[\em (iii)] $G$ is an $\{F_2$, $F_3\}$-free 
                  graph in which every block is a probe complete graph.
\end{itemize}
\end{Theorem}
\begin{figure}[H]
\begin{center}


\begin{tikzpicture}[scale=.55]
\tikzstyle{vP}=[circle,inner sep=1.5pt,fill=black];
\node[vP] (1) at (0,0) {};
\node[vP] (2) at (0,1.5) {};
\node[vP] (3) at (1,2.5) {};
\node[vP] (4) at (2,1.5) {};
\node[vP] (5) at (2,0) {};

\draw (1) -- (2) -- (3) -- (4) -- (5) -- (1) -- (4) -- (2) -- (5);

\node[] (F1) at (1,0) [label=below:$F_1$] {};
\end{tikzpicture}
\qquad\quad
\begin{tikzpicture}[scale=.55]
\tikzstyle{vP}=[circle,inner sep=1.5pt,fill=black];
\node[vP] (1) at (0,1) {};
\node[vP] (2) at (1,0) {};
\node[vP] (3) at (1,2) {};
\node[vP] (4) at (2,1) {};
\node[vP] (5) at (3,0) {};
\node[vP] (6) at (3,2) {};
\node[vP] (7) at (4,1) {};

\draw (1) -- (2) -- (3) -- (4) -- (5) -- (7) -- (6) -- (4) -- (2);
\draw (1) -- (3); \draw (4) -- (7);

\node[] (F2) at (2,0) [label=below:$F_2$] {};
\end{tikzpicture}
\qquad\quad
\begin{tikzpicture}[scale=.55]
\tikzstyle{vP}=[circle,inner sep=1.5pt,fill=black];
\node[vP] (1) at (0,1) {};
\node[vP] (2) at (1,0) {};
\node[vP] (3) at (1,2) {};
\node[vP] (4) at (2,1) {};
\node[vP] (5) at (3,1) {};
\node[vP] (6) at (4,0) {};
\node[vP] (7) at (4,2) {};
\node[vP] (8) at (5,1) {};

\draw (1) -- (2) -- (3) -- (4) -- (5) -- (6) -- (7) -- (8) -- (6);
\draw (1) -- (3); \draw (2) -- (4); \draw (5) -- (7);

\node[] (F3) at (2.5,0) [label=below:$F_3$] {};
\end{tikzpicture}
\end{center}
\caption{Forbidden induced subgraphs for unpartitioned probe block graphs.}
\label{fig:F1234}
\end{figure}
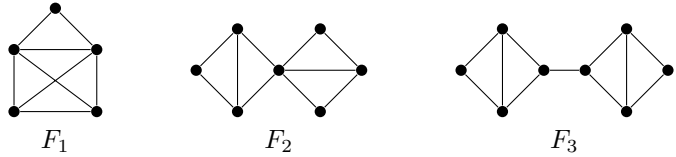
\begin{Theorem}[\cite{LePeng12}]
Partitioned and unpartitioned probe block graphs can be recognized in linear time.
\end{Theorem}

Condition (ii) in Theorem~\ref{thm:probeblock} above can be interpreted as follows: The block structure is described by forbidding $F_1$, the gluing conditions for the blocks are given by forbidding $F_2$ and $F_3$.

\section{$2$-probe complete graphs}\label{sec:2-probecomplete}
For fixed $k$, it was shown in~\cite{ChaHunKloPen} that, by a non-constructive proof, $k$-probe complete graphs can be characterized by at most $2^{k+1}+1$ obstructions, and that $k$-probe complete graphs can be recognized in cubic time. Here we give the complete list of five obstructions for $2$-probe complete graphs, and provide a second characterization of $2$-probe complete graphs. These results imply a linear time recognition algorithm for $2$-probe complete graphs, and are important when discussing $2$-probe block graphs later.

A graph $G=(V,E)$ is called a $(K,X,Y,Z)$-graph, written $G=(K,X,Y,Z)$, if $V$ can be partitioned into disjoint (possibly empty) subsets $K, X, Y, Z$ such that
\begin{itemize}
 \item $K$ is the set of all universal vertices of $G$ (hence, $K$ is a clique),
 \item $X\cup Z$ and $Y\cup Z$ are independent sets,
 \item every vertex in $X$ is adjacent to every vertex in $Y$.
\end{itemize}

Note that if $X=\emptyset$ or $Y=\emptyset$, or ($Z=\emptyset$ and ($|X|\le 1$ or $|Y|\le 1$)), then a $(K,X,Y,Z)$-graph is a complete split graph. Note also that $G$ is a $(K,X,Y,Z)$-graph if and only if the graph obtained from $G$ by deleting all universal vertices has at most one nontrivial connected component which is a complete bipartite graph. Hence, $(K,X,Y,Z)$-graphs can be recognized in linear time.

\begin{Theorem}\label{thm:2-probecomplete}
The following statements are equivalent for any graph $G$:
\begin{itemize}
  \item[\em (i)] $G$ is a $2$-probe complete graph;
  \item[\em (ii)]  $G$ is $\{P_4, 2K_2, K_3+ K_1,(K_2+K_1)\star 2K_1, 2K_1\star 2K_1\star 2K_1\}$-free; 
  \item[\em (iii)] $G$ is a $(K,X,Y,Z)$-graph.
\end{itemize}
\end{Theorem}
\proof

\noindent
(i) $\Rightarrow$ (ii): By inspection, none of $P_4$, $2K_2$, $K_3+ K_1$, $(K_2+K_1)\star 2K_1$ and $2K_1\star 2K_1\star 2K_1$ 
is a $2$-probe complete graph.

\smallskip\noindent
(ii) $\Rightarrow$ (iii): By induction. Let $G$ satisfy (ii). Suppose first, that $G$ is connected. Then, as $G$ is a cograph, $G$ is the join of two smaller graphs, say $G=G_1\star G_2$. By induction, $G_i$ is a $(K_i,X_i,Y_i,Z_i)$-graph, $i=1,2$. If $G_1$ or $G_2$ is a clique, say $G_2$, then $G$ is a $(K_1\cup V(G_2),X_1,Y_1,Z_1)$-graph. If both $G_1$ and $G_2$ are not cliques, then both $G_1$ and $G_2$ are $\{K_2+K_1, C_4\}$-free (otherwise $G$ would have an induced $(K_2+K_1)\star 2K_1$ or $2K_1\star 2K_1\star 2K_1$. 
By Proposition~\ref{prop:completesplit}, for each $i=1,2$, $G_i=(Q_i,S_i)$ is a complete split graph. Hence $G$ is a $(K,X,Y,\emptyset)$-graph with $K=Q_1\cup Q_2$, $X=S_1, Y=S_2$. Suppose now that $G$ is disconnected. We may assume that $G$ is not edgeless. Then $G$ has exactly one nontrivial connected component (as $G$ is $2K_2$-free), say $H$. Note that $H$ is $K_3$-free (as $G$ is $(K_3+K_1)$-free), and hence $H$ is complete bipartite (as $H$ is a connected cograph). Thus $G$ is a $(\emptyset,X,Y,Z)$-graph with $(X,Y)$ being the bipartition of $H$ and $Z=V(G)\setminus V(H)$.

\smallskip\noindent
(iii) $\Rightarrow$ (i): This is obvious by setting $\N_1=X\cup Z$, $\N_2=Y\cup Z$.
\qed

\begin{Corollary}\label{cor:p2-probecomplete}
Let $G=(V,E)$ be a graph with two independent sets $\N_1,\N_2$ and $\Probe=V\setminus\big(\N_1\cup\N_2\big)$.
Then $G=(\Probe,\N_1,\N_2, E)$ is a partitioned $2$-probe complete graph if and only if $G$ is a $(K,X,Y,Z)$-graph such that
$\N_1=X\cup Z$ and $\N_2=Y\cup Z$.
\end{Corollary}

As $(K,X,Y,Z)$-graphs can be recognized in linear time, we obtain:
\begin{Corollary}\label{cor:recog-p2-probecomplete}
Unpartitioned and partitioned $2$-probe complete graphs can be recognized in linear time.
\end{Corollary}

\begin{Corollary}\label{cor:p2-probeblock}
$2$-probe block graphs are distance-hereditary.
\end{Corollary}
\proof\, A slightly stronger statement holds. Each block of a $2$-probe block graph is clearly a $2$-probe complete graph. By Theorem~\ref{thm:2-probecomplete}, each block of a $2$-probe block graph is therefore a cograph.\qed

\section{Partitioned $2$-probe block graphs}\label{sec:p2-probeblock}
Let $G=(V,E)$ be a graph with two given independent set $\N_1, \N_2\subseteq V$. Suppose there exists a set
$E'\subseteq \binom{\N_1}{2}\cup \binom{\N_2}{2}$ such that the graph $G=(V,E\cup E')$ is a block graph, that is, $G$ is a
partitioned $2$-probe block graph with respect to the given independent sets $\N_1, \N_2$.
Then, clearly, the two non-adjacent vertices $x,y$ of every induced diamond in $G$ must belong to one of $\N_1,\N_2$ and $\{x,y\}$ must belong to $E'$. Similarly, any two non-adjacent vertices $x,y$ of every induced $4$-cycle in $G$ must belong to one of $\N_1,\N_2$ and $\{x,y\}$ must belong to $E'$.

In what follows, given a graph $G=(V,E)$ together with two given independent sets $\N_1, \N_2$,
the \emph{enhanced graph} $G^*=(V,E^*)$ is obtained from $G$ by adding all edges between two vertices both in $\N_1$ or both in $\N_2$ that are two vertices of an induced diamond or of an induced $C_4$ in $G$.

Partitioned probe block graphs can be characterized as follows.

\begin{Theorem}\label{thm:p2-probeblock}
Let $G=(V,E)$ be a graph with two independent sets $\N_1,\N_2$ and $\Probe=V\setminus\big(\N_1\cup\N_2\big)$.
Then the following statements are equivalent:
\begin{itemize}
 \item[\em (i)]  $G=(\Probe,\N_1,\N_2, E)$ is a partitioned $2$-probe block graph;
 \item[\em (ii)] Every block $B$ of $G$ is a $(K_B,X_B,Y_B,Z_B)$-graph such that $\N_1\cap B=X_B\cup Z_B$ and $\N_2\cap B=Y_B\cup Z_B$;
 \item[\em (iii)] The enhanced graph $G^*$ of $G$ is a block graph.
\end{itemize}
\end{Theorem}
\proof

\smallskip\noindent
(i) $\Rightarrow$ (ii): Since every block $B$ of $G$ is contained in a block of any block graph embedding of $G$, $B$ is a partitioned $2$-probe complete graph. Hence (ii) follows by Corollary~\ref{cor:p2-probecomplete}.

\smallskip\noindent
(ii) $\Rightarrow$ (iii): Let $G=(V,E)$ satisfy (ii). Then, for every two vertices $x, y$ of $G$, we have the following fact:
\begin{align*}
&\text{Both $x$ and $y$ are in $\N_1$ or in $\N_2$ and belong to an induced diamond or to}\\
&\text{an induced $C_4$ if and only if $x,y$ are non-adjacent vertices in a block of $G$.}
\end{align*}
To see this, note first that one direction is obvious: every diamond and every $C_4$ is contained in a block of $G$. Conversely, let $x, y$ be two non-adjacent vertices in a block $B$ of $G$. As $G$ satisfies (ii), $B=(K_B, X_B, Y_B, Z_B)$ and hence $x,y\in X_B\cup Z_B\subseteq \N_1$ or  $x,y\in Y_B\cup Z_B\subseteq \N_2$. Moreover, as $B$ is 2-connected, it is easy to see that $x, y$ are contained in a diamond or a $C_4$ in $B$. 

Thus, by definition of $G^*=(V,E^*)$, $xy\in E^*\setminus E$ if and only if $x,y$ are non-adjacent vertices of a block in $G$. Therefore, each block of $G^*$ is a clique, that is, $G^*$ is a block graph.

\smallskip\noindent
(iii) $\Rightarrow$ (i): This implication is obvious.
\qed

\medskip
Since the blocks of a graph can be computed in linear time, and $(K,X,Y,Z)$-graphs can be recognized in linear time, Theorem~\ref{thm:p2-probeblock} (ii) implies:
\begin{Corollary}\label{cor:recog-p2-probeblock}
Partitioned $2$-probe block graphs can be recognized in linear time.
\end{Corollary}

\section{Unpartitioned $2$-probe block graphs}\label{sec:2-probeblock}
\begin{figure}[htb]%
\begin{center}
\begin{tikzpicture}[scale=.55]
\tikzstyle{vP}=[circle,inner sep=1.5pt,fill=black];
\node[vP] (1) at (0,1) {};
\node[vP] (2) at (1,0) {};
\node[vP] (34) at (2,1.5) {};
\node[vP] (5) at (1,3) {};
\node[vP] (6) at (0,2) {};

\draw (1) -- (2) -- (34) -- (5) -- (6) -- (1); \draw (1) -- (5); \draw (2) -- (6);

\node[] (B1) at (1,0) [label=below:$B_{1}$] {}; 
\end{tikzpicture}
\quad\,
\begin{tikzpicture}[scale=.55]
\tikzstyle{vP}=[circle,inner sep=1.5pt,fill=black];
\node[vP] (1) at (0,1) {};
\node[vP] (2) at (1,0) {};
\node[vP] (3) at (2,1) {};
\node[vP] (4) at (2,2) {};
\node[vP] (5) at (1,3) {};
\node[vP] (6) at (0,2) {};

\draw (1) -- (2) -- (3) -- (4) -- (5) -- (6) -- (1) -- (4); \draw (3) -- (6);

\node[] (B2) at (1,0) [label=below:$B_{2}$] {}; 
\end{tikzpicture}
\quad\,
\begin{tikzpicture}[scale=.55]
\tikzstyle{vP}=[circle,inner sep=1.5pt,fill=black];
\node[vP] (1) at (0,1) {};
\node[vP] (2) at (1,0) {};
\node[vP] (3) at (2,1) {};
\node[vP] (4) at (2,2) {};
\node[vP] (5) at (1,3) {};
\node[vP] (6) at (0,2) {};

\draw (1) -- (2) -- (3) -- (4) -- (5) -- (6) -- (1) -- (4); \draw (3) -- (6);
\draw (1) -- (5) -- (3);
\draw (4) -- (2) -- (6);

\node[] (B3) at (1,0) [label=below:$B_{3}$] {}; 
\end{tikzpicture}
\quad\,
\begin{tikzpicture}[scale=.55]
\tikzstyle{vP}=[circle,inner sep=1.5pt,fill=black];
\node[vP] (1) at (0,1) {};
\node[vP] (2) at (1,0) {};
\node[vP] (3) at (2,1) {};
\node[vP] (4) at (2,2) {};
\node[vP] (5) at (1,3) {};
\node[vP] (6) at (0,2) {};

\draw (1) -- (2) -- (3) -- (4) -- (5) -- (6) -- (1) -- (3) -- (6) -- (4) -- (1);

\node[] (B4) at (1,0) [label=below:$B_{4}$] {}; 
\end{tikzpicture}
\quad\,
\begin{tikzpicture}[scale=.55]
\tikzstyle{vP}=[circle,inner sep=1.5pt,fill=black];
\node[vP] (1) at (0,1) {};
\node[vP] (2) at (1,0) {};
\node[vP] (3) at (2,1) {};
\node[vP] (4) at (2,2) {};
\node[vP] (5) at (1,3) {};
\node[vP] (6) at (0,2) {};

\draw (1) -- (2) -- (3) -- (4) -- (5) -- (6) -- (1) -- (3) -- (6) -- (4); 
\draw (3) -- (5); \draw (2) -- (6);

\node[] (B5) at (1,0) [label=below:$B_{5}$] {}; 
\end{tikzpicture}
\quad\,
\begin{tikzpicture}[scale=.55]
\tikzstyle{vP}=[circle,inner sep=1.5pt,fill=black];
\node[vP] (1) at (0,1) {};
\node[vP] (2) at (1,0) {};
\node[vP] (3) at (2,1) {};
\node[vP] (4) at (2,2) {};
\node[vP] (5) at (1,3) {};
\node[vP] (6) at (0,2) {};

\draw (1) -- (2) -- (3) -- (4) -- (5) -- (6) -- (1) -- (3) -- (6) -- (4) -- (1);
\draw (2) -- (4); \draw (2) -- (6);

\node[] (B6) at (1,0) [label=below:$B_{6}$] {}; 
\end{tikzpicture}
\end{center}
\caption{Block structure: $2$-connected forbidden induced subgraphs for unpartitioned $2$-probe block graphs.}
\label{fig:B1-6}
\end{figure}
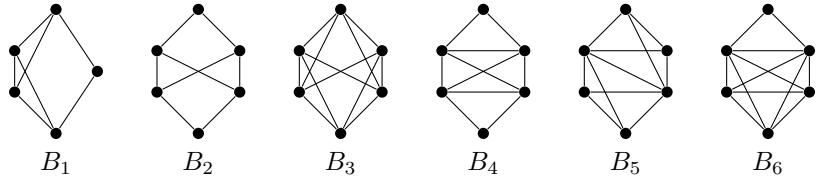

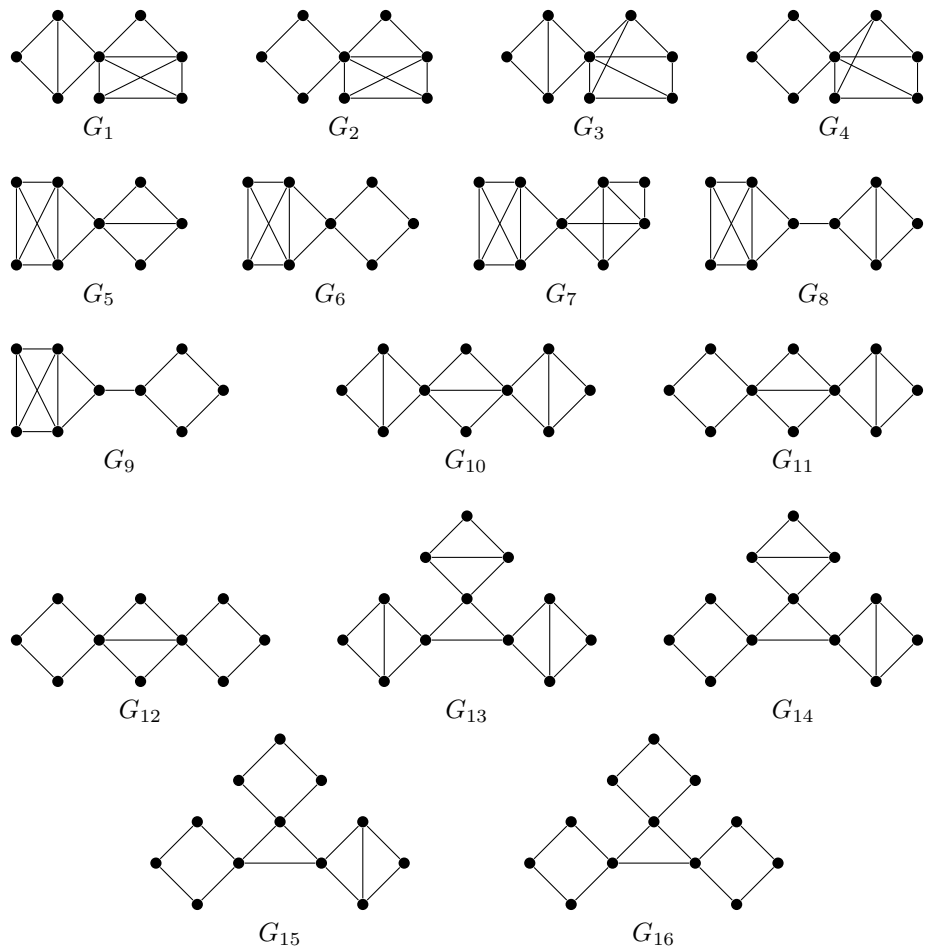
\begin{figure}[H]
\begin{center}
\begin{tikzpicture}[scale=.55]
\tikzstyle{vP}=[circle,inner sep=1.5pt,fill=black];
\node[vP] (1) at (0,1) {};
\node[vP] (2) at (1,0) {};
\node[vP] (3) at (1,2) {};
\node[vP] (4) at (2,1) {};
\node[vP] (5a) at (2,0) {};
\node[vP] (5b) at (4,0) {};
\node[vP] (6) at (3,2) {};
\node[vP] (7) at (4,1) {};

\draw (2) -- (4) -- (3) -- (2) -- (1) -- (3);
\draw (5a) -- (4) -- (6) -- (7) -- (5b) -- (5a) -- (7) -- (4) -- (5b);

\node[] (G1) at (2,0) [label=below:$G_1$] {}; 
\end{tikzpicture}
\hfill
\begin{tikzpicture}[scale=.55]
\tikzstyle{vP}=[circle,inner sep=1.5pt,fill=black];
\node[vP] (1) at (0,1) {};
\node[vP] (2) at (1,0) {};
\node[vP] (3) at (1,2) {};
\node[vP] (4) at (2,1) {};
\node[vP] (5a) at (2,0) {};
\node[vP] (5b) at (4,0) {};
\node[vP] (6) at (3,2) {};
\node[vP] (7) at (4,1) {};

\draw (2) -- (4) -- (3) -- (1) -- (2);
\draw (5a) -- (4) -- (6) -- (7) -- (5b) -- (5a) -- (7) -- (4) -- (5b);

\node[] (G2) at (2,0) [label=below:$G_2$] {}; 
\end{tikzpicture}
\hfill
\begin{tikzpicture}[scale=.55]
\tikzstyle{vP}=[circle,inner sep=1.5pt,fill=black];
\node[vP] (1) at (0,1) {};
\node[vP] (2) at (1,0) {};
\node[vP] (3) at (1,2) {};
\node[vP] (4) at (2,1) {};
\node[vP] (5a) at (2,0) {};
\node[vP] (5b) at (4,0) {};
\node[vP] (6) at (3,2) {};
\node[vP] (7) at (4,1) {};

\draw (2) -- (4) -- (3) -- (2) -- (1) -- (3);
\draw (5a) -- (4) -- (6) -- (7) -- (5b) -- (5a) -- (6); \draw (7) -- (4) -- (5b);

\node[] (G3) at (2,0) [label=below:$G_3$] {}; 
\end{tikzpicture}
\hfill
\begin{tikzpicture}[scale=.55]
\tikzstyle{vP}=[circle,inner sep=1.5pt,fill=black];
\node[vP] (1) at (0,1) {};
\node[vP] (2) at (1,0) {};
\node[vP] (3) at (1,2) {};
\node[vP] (4) at (2,1) {};
\node[vP] (5a) at (2,0) {};
\node[vP] (5b) at (4,0) {};
\node[vP] (6) at (3,2) {};
\node[vP] (7) at (4,1) {};

\draw (2) -- (4) -- (3) -- (1) -- (2);
\draw (5a) -- (4) -- (6) -- (7) -- (5b) -- (5a) -- (6); \draw (7) -- (4) -- (5b);

\node[] (G4) at (2,0) [label=below:$G_4$] {}; 
\end{tikzpicture}
\end{center}
\begin{center}
\begin{tikzpicture}[scale=.55]
\tikzstyle{vP}=[circle,inner sep=1.5pt,fill=black];
\node[vP] (1a) at (0,0) {};
\node[vP] (1b) at (0,2) {};
\node[vP] (2) at (1,0) {};
\node[vP] (3) at (1,2) {};
\node[vP] (4) at (2,1) {};
\node[vP] (5) at (3,0) {};
\node[vP] (6) at (3,2) {};
\node[vP] (7) at (4,1) {};

\draw (1a) -- (2) -- (3) -- (4) -- (5) -- (7) -- (6) -- (4) -- (2);
\draw (1a) -- (3); \draw (4) -- (7); \draw (1a) -- (1b) -- (2); \draw (1b) -- (3);

\node[] (G5) at (2,0) [label=below:$G_5$] {}; 
\end{tikzpicture}
\hfill
\begin{tikzpicture}[scale=.55]
\tikzstyle{vP}=[circle,inner sep=1.5pt,fill=black];
\node[vP] (1a) at (0,0) {};
\node[vP] (1b) at (0,2) {};
\node[vP] (2) at (1,0) {};
\node[vP] (3) at (1,2) {};
\node[vP] (4) at (2,1) {};
\node[vP] (5) at (3,0) {};
\node[vP] (6) at (3,2) {};
\node[vP] (7) at (4,1) {};

\draw (1a) -- (2) -- (3) -- (4) -- (5) -- (7) -- (6) -- (4) -- (2);
\draw (1a) -- (3); 
\draw (1a) -- (1b) -- (2); \draw (1b) -- (3);

\node[] (G6) at (2,0) [label=below:$G_6$] {}; 
\end{tikzpicture}
\hfill
\begin{tikzpicture}[scale=.55]
\tikzstyle{vP}=[circle,inner sep=1.5pt,fill=black];
\node[vP] (1a) at (0,0) {};
\node[vP] (1b) at (0,2) {};
\node[vP] (2) at (1,0) {};
\node[vP] (3) at (1,2) {};
\node[vP] (4) at (2,1) {};
\node[vP] (5) at (3,0) {};
\node[vP] (6) at (3,2) {};
\node[vP] (7) at (4,1) {};
\node[vP] (8) at (4,2) {};

\draw (1a) -- (2) -- (3) -- (4) -- (5) -- (7) -- (6) -- (4) -- (2);
\draw (1a) -- (3); \draw (4) -- (7); \draw (1a) -- (1b) -- (2); \draw (1b) -- (3); 
\draw (5) -- (6) -- (8) -- (7);

\node[] (G7) at (2,0) [label=below:$G_7$] {}; 
\end{tikzpicture}
\hfill
\begin{tikzpicture}[scale=.55]
\tikzstyle{vP}=[circle,inner sep=1.5pt,fill=black];
\node[vP] (1a) at (0,0) {};
\node[vP] (1b) at (0,2) {};
\node[vP] (2) at (1,0) {};
\node[vP] (3) at (1,2) {};
\node[vP] (4) at (2,1) {};
\node[vP] (5) at (3,1) {};
\node[vP] (6) at (4,0) {};
\node[vP] (7) at (4,2) {};
\node[vP] (8) at (5,1) {};

\draw (1a) -- (2) -- (3) -- (4) -- (5) --(6) -- (8) -- (7) -- (5); \draw (6) -- (7);
\draw (1a) -- (3); \draw (2) -- (4); \draw (1a) -- (1b) -- (2); \draw (1b) -- (3);

\node[] (G8) at (2.5,0) [label=below:$G_8$] {}; 
\end{tikzpicture}
\end{center}
\begin{center}
\begin{tikzpicture}[scale=.55]
\tikzstyle{vP}=[circle,inner sep=1.5pt,fill=black];
\node[vP] (1a) at (0,0) {};
\node[vP] (1b) at (0,2) {};
\node[vP] (2) at (1,0) {};
\node[vP] (3) at (1,2) {};
\node[vP] (4) at (2,1) {};
\node[vP] (5) at (3,1) {};
\node[vP] (6) at (4,0) {};
\node[vP] (7) at (4,2) {};
\node[vP] (8) at (5,1) {};

\draw (1a) -- (2) -- (3) -- (4) -- (5) --(6) -- (8) -- (7) -- (5); 
\draw (1a) -- (3); \draw (2) -- (4); \draw (1a) -- (1b) -- (2); \draw (1b) -- (3);

\node[] (G9) at (2.5,0) [label=below:$G_9$] {}; 
\end{tikzpicture}
\hfill\quad\,
\begin{tikzpicture}[scale=.55]
\tikzstyle{vP}=[circle,inner sep=1.5pt,fill=black];
\node[vP] (1) at (0,1) {};
\node[vP] (2) at (1,0) {};
\node[vP] (3) at (1,2) {};
\node[vP] (4) at (2,1) {};
\node[vP] (5) at (3,0) {};
\node[vP] (6) at (3,2) {};
\node[vP] (7) at (4,1) {};
\node[vP] (8) at (5,0) {};
\node[vP] (9) at (5,2) {};
\node[vP] (10) at (6,1) {};

\draw (2) -- (4) -- (3) -- (1) -- (2) -- (3);
\draw (4) -- (5) -- (7) -- (6) -- (4) -- (7);
\draw (8) -- (10) -- (9) -- (7) -- (8) -- (9);

\node[] (G12) at (3,0) [label=below:$G_{10}$] {}; 
\end{tikzpicture}
\hfill
\begin{tikzpicture}[scale=.55]
\tikzstyle{vP}=[circle,inner sep=1.5pt,fill=black];
\node[vP] (1) at (0,1) {};
\node[vP] (2) at (1,0) {};
\node[vP] (3) at (1,2) {};
\node[vP] (4) at (2,1) {};
\node[vP] (5) at (3,0) {};
\node[vP] (6) at (3,2) {};
\node[vP] (7) at (4,1) {};
\node[vP] (8) at (5,0) {};
\node[vP] (9) at (5,2) {};
\node[vP] (10) at (6,1) {};

\draw (2) -- (4) -- (3) -- (1) -- (2); 
\draw (4) -- (5) -- (7) -- (6) -- (4) -- (7);
\draw (8) -- (10) -- (9) -- (7) -- (8) -- (9);

\node[] (G13) at (3,0) [label=below:$G_{11}$] {}; 
\end{tikzpicture}
\end{center}
\begin{center}
\begin{tikzpicture}[scale=.55]
\tikzstyle{vP}=[circle,inner sep=1.5pt,fill=black];
\node[vP] (1) at (0,1) {};
\node[vP] (2) at (1,0) {};
\node[vP] (3) at (1,2) {};
\node[vP] (4) at (2,1) {};
\node[vP] (5) at (3,0) {};
\node[vP] (6) at (3,2) {};
\node[vP] (7) at (4,1) {};
\node[vP] (8) at (5,0) {};
\node[vP] (9) at (5,2) {};
\node[vP] (10) at (6,1) {};

\draw (2) -- (4) -- (3) -- (1) -- (2); 
\draw (4) -- (5) -- (7) -- (6) -- (4) -- (7);
\draw (8) -- (10) -- (9) -- (7) -- (8); 

\node[] (G14) at (3,0) [label=below:$G_{12}$] {}; 
\end{tikzpicture}
\hfill
\begin{tikzpicture}[scale=.55]
\tikzstyle{vP}=[circle,inner sep=1.5pt,fill=black];
\node[vP] (1) at (0,1) {};
\node[vP] (2) at (1,0) {};
\node[vP] (3) at (1,2) {};
\node[vP] (4) at (2,1) {};
\node[vP] (5) at (3,2) {};
\node[vP] (6) at (2,3) {};
\node[vP] (7) at (3,4) {};
\node[vP] (8) at (4,3) {};
\node[vP] (9) at (4,1) {};
\node[vP] (10) at (5,0) {};
\node[vP] (11) at (5,2) {};
\node[vP] (12) at (6,1) {};

\draw (2) -- (4) -- (3) -- (1) -- (2) -- (3);
\draw (5) -- (4) -- (9) -- (5);
\draw (6) -- (7) -- (8) -- (5) -- (6) -- (8);
\draw (11) -- (12) -- (10) -- (9) -- (11) -- (10);

\node[] (G15) at (3,0) [label=below:$G_{13}$] {}; 
\end{tikzpicture}
\hfill
\begin{tikzpicture}[scale=.55]
\tikzstyle{vP}=[circle,inner sep=1.5pt,fill=black];
\node[vP] (1) at (0,1) {};
\node[vP] (2) at (1,0) {};
\node[vP] (3) at (1,2) {};
\node[vP] (4) at (2,1) {};
\node[vP] (5) at (3,2) {};
\node[vP] (6) at (2,3) {};
\node[vP] (7) at (3,4) {};
\node[vP] (8) at (4,3) {};
\node[vP] (9) at (4,1) {};
\node[vP] (10) at (5,0) {};
\node[vP] (11) at (5,2) {};
\node[vP] (12) at (6,1) {};

\draw (2) -- (4) -- (3) -- (1) -- (2); 
\draw (5) -- (4) -- (9) -- (5);
\draw (6) -- (7) -- (8) -- (5) -- (6) -- (8);
\draw (11) -- (12) -- (10) -- (9) -- (11) -- (10);

\node[] (G16) at (3,0) [label=below:$G_{14}$] {}; 
\end{tikzpicture}
\end{center}
\vspace*{-2em}
\begin{center}
\begin{tikzpicture}[scale=.55]
\tikzstyle{vP}=[circle,inner sep=1.5pt,fill=black];
\node[vP] (1) at (0,1) {};
\node[vP] (2) at (1,0) {};
\node[vP] (3) at (1,2) {};
\node[vP] (4) at (2,1) {};
\node[vP] (5) at (3,2) {};
\node[vP] (6) at (2,3) {};
\node[vP] (7) at (3,4) {};
\node[vP] (8) at (4,3) {};
\node[vP] (9) at (4,1) {};
\node[vP] (10) at (5,0) {};
\node[vP] (11) at (5,2) {};
\node[vP] (12) at (6,1) {};

\draw (2) -- (4) -- (3) -- (1) -- (2); 
\draw (5) -- (4) -- (9) -- (5);
\draw (6) -- (7) -- (8) -- (5) -- (6); 
\draw (11) -- (12) -- (10) -- (9) -- (11) -- (10);

\node[] (G17) at (3,0) [label=below:$G_{15}$] {}; 
\end{tikzpicture}
\qquad\qquad
\begin{tikzpicture}[scale=.55]
\tikzstyle{vP}=[circle,inner sep=1.5pt,fill=black];
\node[vP] (1) at (0,1) {};
\node[vP] (2) at (1,0) {};
\node[vP] (3) at (1,2) {};
\node[vP] (4) at (2,1) {};
\node[vP] (5) at (3,2) {};
\node[vP] (6) at (2,3) {};
\node[vP] (7) at (3,4) {};
\node[vP] (8) at (4,3) {};
\node[vP] (9) at (4,1) {};
\node[vP] (10) at (5,0) {};
\node[vP] (11) at (5,2) {};
\node[vP] (12) at (6,1) {};

\draw (2) -- (4) -- (3) -- (1) -- (2); 
\draw (5) -- (4) -- (9) -- (5);
\draw (6) -- (7) -- (8) -- (5) -- (6); 
\draw (11) -- (12) -- (10) -- (9) -- (11); 

\node[] (G18) at (3,0) [label=below:$G_{16}$] {}; 
\end{tikzpicture}
\end{center}
\caption{Gluing conditions: Forbidden induced subgraphs for unpartitioned $2$-probe block graphs.} 
\label{fig:G1-16}
\end{figure}

In this section, we characterize $2$-probe block graphs in terms of their block structure and gluing conditions. The characterization reminds the one of $1$-probe block graphs (Theorem~\ref{thm:probeblock}), but it is considerably more involved. It turns out that the blocks are $2$-probe complete graphs and can be described by six forbidden induced subgraphs depicted in Figure~\ref{fig:B1-6}, and the gluing conditions can be expressed in terms of the other sixteen forbidden induced subgraphs depicted in Figure~\ref{fig:G1-16}.

\begin{Theorem}\label{thm:2-probeblock}
The following statements are equivalent for any graph $G$:
\begin{itemize}
\item[\em (i)] $G$ is a $2$-probe block graph;
\item[\em (ii)] $G$ is a $\{B_1, \ldots, B_6, G_1,\ldots, G_{16}\}$-free distance-hereditary graph;
\item[\em (iii)] $G$ is a $\{G_1, \ldots, G_{16}\}$-free 
                 graph in which every block is a $2$-probe complete graph.
\end{itemize}
\end{Theorem}

\proof

\smallskip\noindent
(i) $\Rightarrow$ (ii): By Corollary~\ref{cor:p2-probeblock}, $G$ is distance-hereditary. By inspection, none of $B_1$, \ldots, $B_6$, $G_1$, \ldots, $G_{16}$ is a $2$-probe block graph.

\smallskip\noindent
(ii) $\Rightarrow$ (iii): Let $G$ satisfy (ii), and let $B$ be a block of $G$. Then $B$ is a $2$-connected distance-hereditary graph without induced $B_1$, \ldots, $B_6$. We claim that $B$ is $\{P_4, 2K_2, K_3+K_1\}$-free.

Assume first that $B$ contains an induced subgraph $P_4$, say $abcd$. Since $B$ is $2$-connected distance-hereditary, there is a vertex $x$ adjacent to $a$ and $c$ (hence non-adjacent to $d$), and another vertex $y$ adjacent to $b$ and $d$ (hence non-adjacent to $a$). Let $H$ be the subgraph of $B$ induced by $a,b,c,d,x$ and $y$. As $H$ is not a domino, one of the edges $xy, xb, yc$ must exist. Then, an easy case analysis shows that $H$ is isomorphic to $B_2$ or to $B_4$, or $H-a$ or $H-d$ is a house or a gem. Thus, $B$ is a cograph. Since $B$ is $2$-connected, $B=H_1\star H_2$. Assume next that $B$ contains an induced subgraph $F\in\{2K_2,K_3+K_1\}$. Then $F$ is contained in $H_1$, say, and $|V(H_2)|=1$ (otherwise there would be a $B_1$ or $B_5$ in case $F=2K_2$, or a $B_1$ or $B_6$ in case $F=K_3+K_1$). Therefore $H_1$ is connected and hence $H_1=H_{11}\star H_{12}$ with $F$ being contained in $H_{11}$, say. But now there is a $B_5$ or a $B_6$ induced by $F$, $H_2$ and a vertex in $H_{12}$.

Thus, $B$ is $\{P_4, 2K_2, K_3+K_1\}$-free, as claimed. Note that $B_1=(K_2+K_1)\star 2K_1$ and $B_3=2K_1\star 2K_1\star 2K_1$. 
Hence, by Theorem~\ref{thm:2-probecomplete}, $B$ is a $2$-probe complete graph.

\smallskip\noindent
(iii) $\Rightarrow$ (i): We prove a slightly stronger claim that every graph $H$ satisfying (iii) admits two (possibly empty) independent sets $\N_1$ and $\N_2$ such that
\begin{equation}\label{eq1}
\text{$H=(\Probe,\N_1, \N_2,E)$ is a partitioned $2$-probe block graph,}
\end{equation}
\begin{equation}\label{eq2}
\begin{split}
&\text{$\forall\, i=1,2$, $\forall\, v\in\N_i$, there is another vertex $v'\in\N_i$ such that $v$}\\
&\text{and $v'$ are degree-$2$ vertices of an induced $C_4$ or diamond in $H$,}
\end{split}
\end{equation}
\begin{equation}\label{eq3}
\begin{split}
&\text{every vertex $v\in\N_1\cap \N_2$ is the degree-$2$ vertex of some}\\
&\text{induced $F_1$ (see Figure~\ref{fig:F1234}) in $H$,}
\end{split}
\end{equation}
and,
\begin{equation}\label{eq4}
\begin{split}
&\text{$\forall\, v\in\Probe$, $\forall\, x\in N(v)\cap\N_1, \forall\, y\in N(v)\cap \N_2$, it holds that:}\\
&\text{$x\in\N_2$ or $y\in\N_1$ or $xy\in E(H)$.}
\end{split}
\end{equation}

We will prove this claim by induction. Let $G$ satisfy (iii).
If $G$ is itself a block, then by assumption, $G$ is a $2$-probe complete graph, and by Theorem~\ref{thm:2-probecomplete}, $G$ is a $(K,X,Y,Z)$-graph. If $G$ is a clique, set $\N_1=\N_2=\emptyset$. If $X\not=\emptyset$ and $Y\not=\emptyset$, set $\N_1=X\cup Z$, $\N_2=Y\cup Z$. Finally, if $X=\emptyset$ or $Y=\emptyset$, set $\N_1=X\cup Y\cup Z$, $\N_2=\emptyset$. It is clear, by the $2$-connectedness of $G$, the properties (\ref{eq1}), (\ref{eq2}), (\ref{eq3}) and (\ref{eq4}) hold in this case.

So, consider an end-block $B$ of $G$ and let $v$ be the cut-vertex of $G$ in $B$. Let $H=G-(V(B)\setminus\{v\})$. By induction, $H$ admits independent sets $\N_1', \N_2'$ satisfying (\ref{eq1}), (\ref{eq2}), (\ref{eq3}) and (\ref{eq4}).

If $B$ is a clique, then clearly $\N_1:= \N_1'$ and $\N_2:=\N_2'$ are independent sets of $G$ satisfying (\ref{eq1}), (\ref{eq2}), (\ref{eq3}) and (\ref{eq4}), and we are done.

So, we may assume that $B$ is not a clique. By assumption, $B$ is a $2$-probe complete graph. Write $B=(K,X,Y,Z)$, where $K$ is the set of all universal vertices of $B$ and $Z$ is the set of all isolated vertices of $B-K$. Then, as $B$ is $2$-connected and not complete, 
\begin{itemize}
\item $Z\not=\emptyset$ and $|K|\ge 2$, or else 
\item $|X|\ge 2$ and $|Y|\ge 2$. 
\end{itemize}
Moreover, by definition of $Z$ and $K$, 
\begin{itemize}
\item $X=\emptyset$ if and only if $Y=\emptyset$, and 
\item if $Z=\emptyset$, then $|X|\ge 2, |Y|\ge 2$.
\end{itemize}

\smallskip\noindent
\textsc{Case 1.}\, $v\in K$.\, 
In this case, $v$ belongs to an induced diamond $D$ in $B$ with $\deg_D(v)=3$. This can be seen as follows: In case $|X|\ge 2$ and $|Y|\ge 2$, $v$, two vertices in $X$ and a vertex in $Y$ together induced such a diamond $D$. In other case, $Z\not=\emptyset$ and $|K|\ge 2$. Let $w$ be a vertex in $K\setminus\{v\}$. If $|Z|\ge 2$, then $v, w$ and two vertices in $Z$ induce such a diamond $D$. If $|Z|=1$, then $X\not=\emptyset$ (as $B$ is not complete) and $v, w$, the vertex in $Z$ and a vertex in $X$ induce such a diamond $D$.  

Assume first that $v\in\N_1'\cap\N_2'$. By (\ref{eq3}), $v$ is the degree-$2$ vertex of some induced $F_1$ in $H$. Then, as the cut-vertex $v$ is the only common vertex of this $F_1$ and $D$, this $F_1$ and $D$ together induce a $G_5$, a contradiction. Thus, $v\not\in\N_1'\cap\N_2'$.

Assume next that $v\in\N_1'\cup\N_2'$, say $v\in\N_1'$ but $v\not\in \N_2'$. Then, by
(\ref{eq2}), $v$ belongs to a $4$-vertex induced subgraph $H'$ in $H$ that is a $C_4$ or a diamond with $\deg_{H'}(v)=2$. Hence $Z\not=\emptyset$  (otherwise $H'$, two vertices in $X$ and two vertices in $Y$ together would induce a $G_3$ or $G_4$), and $X=\emptyset$ and $Y=\emptyset$ (otherwise $H'$, a vertex in $Z$, a vertex in $K\setminus\{v\}$, a vertex in $X$ and a vertex in $Y$ together would induce a $G_1$ or $G_2$). Then, as $v\not\in \N_2'$, $\N_1:= \N_1'$ and $\N_2:= \N_2'\cup Z$ are independent sets in $G$ satisfying (\ref{eq1}), (\ref{eq2}), (\ref{eq3}) and (\ref{eq4}), and we are done.

Thus, we may assume that $v\not\in\N_1'\cup\N_2'$. In this case, define $\N_1$ and $\N_2$ as follows. If $X\not=\emptyset$ and $Y\not=\emptyset$, then $\N_1:=\N_1'\cup X\cup Z$ and $\N_2:=\N_2'\cup Y\cup Z$. Otherwise (recall that $X=\emptyset$ if and only if $Y=\emptyset$), $\N_1:=\N_1'\cup Z$ and $\N_2:=\N_2'$. Clearly, $\N_1$ and $\N_2$ are independent sets of $G$ satisfying (\ref{eq1}), (\ref{eq2}), (\ref{eq3}) and (\ref{eq4}).

Case 1 is settled.

\medskip\noindent
\textsc{Case 2}\, $v\not\in K$.\, 
In this case, $v$ belongs to a $4$-vertex induced subgraph $F$ in $B$ such that $F$ is a diamond or a $C_4$ and $\deg_F(v)=2$. Moreover, if $Z=\emptyset$, then $F=C_4$. This can be seen as follows: Suppose first $v\in X\cup Y$, say $v\in X$. If $Z=\emptyset$, then $|X|\ge 2, |Y|\ge 2$, and $v$, another vertex in $X$ and two vertices in $Y$ together induce such an $F=C_4$. If $Z\not=\emptyset$, then $|K|\ge 2$, and $v$, two vertices in $K$ and a vertex in $Z$ induce such a diamond $F$. Suppose next $v\in Z$. If $|Z|\ge 2$, then $v$, another vertex in $Z$ and two vertices in $K$ together induce such a diamond $F$. If $Z=\{v\}$, then $X\not=\emptyset$ (as $B$ is not complete), and $v$, two vertices in $K$ and a vertex in $X$ together induce such a diamond $F$.  

We first show that
\begin{equation*}
\text{$v$ has no neighbor in $\N_1'\cap N_2'$.}
\end{equation*}
For otherwise, let $x\in N(v)\cap\N_1'\cap N_2'$. By (\ref{eq3}), $x$ belongs to an induced subgraph $H'=F_1$ in $H$ with $\deg_{H'}(x)=2$. Let $B'$ be the block of $H$ containing $H'$. Now, if $v\not\in B'$, then $F$ and $H'$ together induce a $G_8$ or $G_9$.
So, let $v\in B'$. Then, by (\ref{eq1}), $v$ is a universal vertex in $B'$ (as $v\not\in\N_1'\cup\N_2'$) and hence $v$ is not one of the two degree-$3$ vertices $y_1, y_2$ of $H'$ (which belong to $\N_1'\cup\N_2'$).
Thus, $v\not\in H'$ or $v$ is one of the degree-$4$ vertices of $H'$. Then $F, x, y_1, y_2$ and a degree-$4$ vertex in $H'$ different from $v$ together induce a $G_{1}$ or $G_{2}$. In any case, we have a contradiction, hence $v$ cannot have a neighbor in $\N_1'\cap N_2'$, as claimed.

\medskip
Next we show that
\begin{equation*}
\text{$N(v)\cap \N_1'=\emptyset$ or $N(v)\cap \N_2'=\emptyset$.}
\end{equation*}
For otherwise let $x_1\in N(v)\cap \N_1'$ and $x_2\in N(v)\cap \N_2'$. Since $v$ has no neighbor in $\N_1'\cap \N_2'$, $x_1\not\in\N_2', x_2\not\in\N_1'$ (in particular, $x_1\not=x_2$). Hence, by (\ref{eq4}), $x_1x_2\in E(H)$. By (\ref{eq2}), $x_i$ belongs to a $4$-vertex induced subgraph $H_i$ in $H$ that is a $C_4$ or a diamond with $\deg_{H_i}(x_i)=2$, $i=1,2$, and the vertex $x_i'$ in $H_i$ non-adjacent to $x_i$ also belongs to $\N_i'$. Let $B_i$ be the blocks of $H$ containing $H_i$. Now, if $v\not\in B_1\cup B_2$, then $B_1\cap B_2=\emptyset$ and $F$, $H_1, H_2$ together induce a $G_{13}, G_{14}, G_{15}$ or $G_{16}$.
If $v\in B_1\cap B_2$, then in particular $B_1= B_2$, and by (\ref{eq1}), $v$ is a universal vertex in $B_1=B_2$, hence (as $v$ has no neighbor in $\N_1'\cap \N_2'$), $x_1'\not\in\N_2'$, $x_2'\not\in\N_1'$. Therefore, by (\ref{eq1}), $x_1, x_2, x_1', x_2'$ must induce a $C_4$ in $B_1=B_2$. But then $F$ and $x_1, x_2, x_1', x_2'$ together induce a $G_3$ or $G_4$.
So, let $v\in B_1$ and $v\not\in B_2$, say. Then $B_1\cap B_2=\{x_2\}$.
As before, by (\ref{eq1}), $v$ is a universal vertex in $B_1$, hence $x_1'\not\in\N_2'$. Therefore $x_1'$ must be adjacent to $x_2$. But then $F$, $x_1, x_1'$, and $H_2$ together induce a $G_{10}, G_{11}$ or $G_{12}$. In any case, we have a contradiction, hence $v$ cannot have neighbors in both $\N_1'$ and $\N_2'$, as claimed.

We now distinguish two subcases.

\medskip\noindent
\textsc{Subcase 2.1.}\, $N(v)\cap\N_1'=N(v)\cap \N_2'=\emptyset$.

Assume first that $v\in\N_1'\cap \N_2'$ and $v\not\in Z$. Then $v\in X\cup Y$ and thus both $X$ and $Y$ are nonempty; let $v\in X$, say. By (\ref{eq3}), $v$ is the degree-$2$ vertex of some induced $H'=F_1$ in $H$. Now, if $Z=\emptyset$ then $F$ is a $C_4$ and therefore $F$ and $H'$ together induce a $G_6$. If $Z\not=\emptyset$ then a vertex in $Z$, two vertices in $K$, a vertex in $Y$ and $H'$ together induce a $G_7$. In any case we have a contradiction.

Thus, $v\in Z$ or $v\not\in\N_1'$ or $v\not\in \N_2'$. Now, define $\N_1$ and $\N_2$ as follows.

Suppose $v\in Z$. If $X\not=\emptyset$ and $Y\not=\emptyset$, $\N_1:=\N_1'\cup X\cup Z$ and $\N_2:=\N_2'\cup Y\cup Z$. Otherwise, $\N_1:=\N_1'\cup Z$ and $\N_2:=\N_2'$.

Suppose $v\not\in Z$ and say $v\not\in\N_1'$. If $v\in X$, $\N_1:=\N_1'\cup Y\cup Z$ and $\N_2:=\N_2'\cup X\cup Z$. Otherwise, $\N_1:=\N_1'\cup X\cup Z$ and $\N_2:=\N_2'\cup Y\cup Z$.

Then, as $N(v)\cap\N_1'=N(v)\cap \N_2'=\emptyset$, $\N_1$ and $\N_2$ are independent sets of $G$, and clearly, satisfy (\ref{eq1}), (\ref{eq2}), (\ref{eq3}) and (\ref{eq4}).

\medskip\noindent
\textsc{Subcase 2.2.}\, $N(v)\cap\N_1'=\emptyset$ and $N(v)\cap \N_2'\not=\emptyset$, say. Then $v\not\in\N_2'$.

Assume first that $v\in Z$. Then $X=Y=\emptyset$. Otherwise, consider a vertex $x\in N(v)\cap \N_2'$. By (\ref{eq2}), $x$ belongs to a $4$-vertex induced subgraph $H'$ in $H$ that is a $C_4$ or a diamond with $\deg_{H'}(x)=2$. Let $B'$ be the block of $H$ containing $H'$.
If $v\not\in B'$, then $v$, two vertices in $K$, a vertex in $X$, a vertex in $Y$ and $H'$ together induce a $G_8$ or $G_9$.
If $v\in B'$, then by (\ref{eq1}), $v$ is a universal vertex in $B'$, or $v\in\N_1'$ and $H'$ is a $C_4$. Recall that the vertex $x'$ in $H'$ nonadjacent to $x$ also belongs to $\N_2'$. Let $y$ be a vertex in $H'-\{x,x',v\}$. Note that $v, x,y,x'$ induce a diamond in $B'$ (if $v$ is universal in $B'$) or a $C_4$ (if $v,y\in\N_1'$). Now, two vertices in $K$, a vertex in $X$, a vertex in $Y$ and $v,x,y,x'$ together induce a $G_5$ or $G_6$.

Thus, if $v\in Z$ then $X=Y=\emptyset$, and, as $N(v)\cap\N_1'=\emptyset$, $\N_1:=\N_1'\cup Z$ and $\N_2:=\N_2'$ are independent sets of $G$ satisfying (\ref{eq1}), (\ref{eq2}), (\ref{eq3}) and (\ref{eq4}).

So, we may assume that $v\in X\cup Y$. Then, if $v\in X$, set $\N_1:=\N_1'\cup X\cup Z$ and $\N_2:=\N_2'\cup Y\cup Z$. If $v\in Y$, set $\N_1:=\N_1'\cup Y\cup Z$ and $\N_2:=\N_2'\cup X\cup Z$. It is clear that, as $N(v)\cap\N_1'=\emptyset$, in any case, $\N_1$ and $\N_2$ are independent sets of $G$ satisfying (\ref{eq1}), (\ref{eq2}), (\ref{eq3}) and (\ref{eq4}).

\medskip
Case 2 is settled. The proof of (iii) $\Rightarrow$ (i) is complete, hence Theorem~\ref{thm:2-probeblock}.\qed

\section{A linear time recognition of unpartitioned $2$-probe block graphs}\label{sec:reog}
Based on Theorem~\ref{thm:2-probeblock} and its proof, we describe briefly in this section how to recognize in linear time whether a given graph is a $2$-probe block graph, and if so, to output a partition into probes and non-probes.

\begin{enumerate}
\item\label{step1}
Compute the blocks and the cut-vertices of $G$.
\item\label{step2}
For each non-complete block $B$ of $G$, let $K_B$ be the set of its universal vertices and $Z_B$ be the set of all isolated vertices of $B-K_B$.\\
If $B-(K_B\cup Z_B)$ is not a complete bipartite graph, then output ``NO'' and STOP, meaning that $G$ is not a $2$-probe block graph.\\
Otherwise, let $(X_B,Y_B)$ be the bipartition of the complete bipartite graph; possibly empty.\\
Let $\mathcal{B}$ be the set of all blocks $B=(K_B, X_B,Y_B,Z_B)$ and their cut-vertices of $G$.
\item\label{step3}
Follow the proof of Theorem~\ref{thm:2-probeblock}, part (iii) $\Rightarrow$ (i), to compute the candidates $\N_1$ and $\N_2$.
\item\label{step4}
Check if $\N_1$ and $\N_2$ are independent sets of $G$. If not, output ``NO'' and STOP. Otherwise check if $G=(\Probe,\N_1,\N_2,E)$ is a partitioned $2$-probe block graph. If this is the case, output ``YES'', meaning $G$ is a $2$-probe block graph. If not, output ``NO''.
\end{enumerate}

Given Theorem~\ref{thm:2-probeblock} and its proof, the correctness is clear. It is also clear that each step, except step~\ref{step3}, can be implemented with linear time complexity.
Step~\ref{step3} can be described more precisely by the following procedure \texttt{FindNonprobes}.

\begin{table}[hbt]
  \begin{center}
    Procedure \texttt{FindNonprobes}$(\mathcal{B}; \N_1,\N_2)$\\[.7ex]
    \fbox{\,
       \begin{minipage}{\textwidth}
        \begin{tabbing}
          \OUTPUT\=\kill
          \INPUT \>  set $\mathcal{B}$ of blocks.\\ 
          \OUTPUT\>  sets of vertices $\N_1, \N_2$.\\[1ex]
          (N.)\hspace*{-1.5ex}\=\ei\=\ei\=\ei\=\ei\=\ei\=\ei\=\ei\=\ei\=\ei\=\ei\kill
          1.\>\> \IF $\mathcal{B}=\{B\}$ \THEN\\
          2.\>\>\> if $B$ is a clique, set $\N_1:=\emptyset; \N_2:= \emptyset$.\\
          3.\>\>\> if $X_B\not=\emptyset$, set $\N_1:=X_B\cup Z_B; \N_2:= Y_B\cup Z_B$.\\
          4.\>\>\> if $X_B=Y_B=\emptyset$, set $N_1:= Z_B, \N_2:=\emptyset$.\\
          5.\>\> \ELSE\\
          6.\>\>\> let $B\in\mathcal{B}$ be an end-block with cut-vertex $v$.\\
          7.\>\>\> set $\mathcal{B}':= \mathcal{B}\setminus\{B\}$ and call \texttt{FindNonprobes}$(\mathcal{B}'; \N_1',\N_2')$.\\
          8.\>\>\> \IF $B$ is a clique \THEN\\
          9.\>\>\>\> $\N_1:=\N_1'; \N_2:= \N_2'$\\
         10.\>\>\> \ELSE\\
         11.\>\>\>\>  follow the proof of Theorem~\ref{thm:2-probeblock} (iii) $\Rightarrow$ (i).\\
         12.\>\>\>\>  if $v\in K_B$, compute $\N_1$ and $\N_2$ according to Case 1.\\
         13.\>\>\>\>  if $v\not\in K_B$, compute $\N_1$ and $\N_2$ according to Case 2.\\
         14.\>\>\> \ENDIF\\
         15.\>\> \ENDIF
        \end{tabbing}
      \end{minipage}
    \quad}
  \end{center}
  \label{alg:FindNonprobes}
\end{table}

Let $t(\mathcal{B})$ denote the time needed by \texttt{FindNonprobes}$(\mathcal{B}; \N_1,\N_2)$. Note that each block $B$ has $|E(B)|\ge |V(B)|$, unless $B\in\{K_1,K_2\}$. Thus, we have $t(\{B\})=O(|E(B)|)+O(1)$ and $t(\mathcal{B})=t(\mathcal{B}\setminus\{B\})+O(|E(B)|)+O(1)$. Therefore, $t(\mathcal{B})=\sum_{B\in\mathcal{B}}( O(|E(B)|)+O(1))= O(|E(G)|+|V(G)|)$.
To sum up, we have:
\begin{Theorem}\label{thm:reg-2-probeblock}
Unpartitioned $2$-probe block graphs can be recognized in linear time.
\end{Theorem}

We note that our algorithm is optimal in the following senses: If $G$ is a $2$-probe block graph, the partition into $(\Probe,\N_1,\N_2)$ is minimal, {\em i.e.}, the corresponding block graph embedding has minimal number of new edges. Moreover, if $G$ is a $1$-probe block graph, then the algorithm will output $\N_2=\emptyset$.

\section{Conclusion}\label{sec:conclusion}
With Theorems \ref{thm:2-probecomplete}, \ref{thm:p2-probeblock}, and~\ref{thm:2-probeblock} we have given good characterizations of $2$-probe complete graphs and $2$-probe block graphs. This might be a
first step towards the solution of the challenging problems of characterizing and recognizing $k$-probe complete graphs and $k$-probe block graphs for any $k\ge 3$.

These problems seem be very difficult even for certain restricted graph classes.
For instance, it is not clear which cographs are $k$-complete graphs, given $k$.

\medskip\noindent
\textbf{Problem 1.}
{\em Let $k\ge 3$. Characterize cographs that are $k$-probe complete graphs, respectively, $k$-probe block graphs.
}

\medskip\noindent
Note that any graph $G=(V,E)$ is a $k$-probe complete graph for some $k$; for instance, $k=\binom{|V|}{2} - |E|$.
In~\cite{ChaHunKloPen}, it is proved that determining the smallest integer $k$ such that
$G$ is a $k$-probe complete graph is NP-hard. Indeed, it is equivalent to determine the smallest number of
cliques in $\overline{G}$ that cover the edges of $\overline{G}$, which is a well-known NP-hard problem~\cite{Orlin77}.

We remark that, in connection to Problem~1, it is still unknown how, given a cograph $G$, to compute the smallest integer $k$ such that $G$ is a $k$-probe complete graph. Indeed, Ton Kloks posed the following problem in 2007.

\medskip\noindent
\textbf{Problem 2 (Kloks~\cite{Kloks07,Kloks-cstheory})}
{\em
Given a cograph $G$, determine the smallest integer $k$ such that $G$ is a $k$-probe complete graph.
}

\medskip\noindent
Since cographs are self-complementary, Kloks' problem is equivalent to
determine the minimum number of cliques that cover the edges of a given cograph $G$.



\end{document}